\documentclass{amsart}
\usepackage{amsrefs,amssymb,amsmath,amsthm,amsfonts,graphicx,epsfig}
\usepackage[utf8]{inputenc}
\usepackage{tikz}
\usetikzlibrary{matrix,arrows,decorations.pathmorphing}

\usepackage{hyperref}
\hypersetup{
    colorlinks,
    citecolor=black,
    filecolor=black,
    linkcolor=black,
    urlcolor=black
}

\theoremstyle{plain}
\newtheorem{teo}{Theorem}[section]
\newtheorem{cor}[teo]{Corollary}
\newtheorem{lemma}[teo]{Lemma}
\newtheorem{prop}[teo]{Proposition}

\newtheorem{obs}[teo]{Remark}

\theoremstyle{definition}

\DeclareMathOperator{\dif}{d}
\DeclareMathOperator{\vertex}{\mathcal{V}}
\DeclareMathOperator{\sing}{Fix}

\DeclareMathOperator{\diam}{diam}

\DeclareMathOperator{\dist}{dist}

\newcommand{\Suno}{\mathbb{S}^1}

\newcommand{\Sdos}{\mathbb{S}^2}
\newcommand{\Stres}{\mathbb{S}^3}

\newcommand{\rt}{\mathcal{R}}

\newcommand{\sect}{H}

\newcommand{\R}{\mathbb R}
\newcommand{\D}{\mathbb D}

\newcommand{\Z}{\mathbb Z}

\newcommand{\T}{\mathcal T}
\newcommand{\Sup}{\mathcal S}
\newcommand{\radio}{\tau}

\renewcommand{\epsilon}{\varepsilon}
\author{Alfonso Artigue}
\address{Departamento de Matem\'atica y Estad\'\i stica del Litoral, 
Universidad de la Rep\'ublica, Gral. Rivera 1350, Salto, Uruguay.}

\title{Expansive Flows of the Three-Sphere}

\begin{document}
\maketitle
\begin{abstract} 
In this article we show that the three-dimensional sphere admits {transitive} expansive flows in the sense of Komuro with hyperbolic 
equilibrium points. 
The result is based on a construction that allows us to see the 
geodesic flow of a hyperbolic three-punctured two-dimensional sphere
as the flow of a smooth vector field on the three-dimensional sphere.
\end{abstract}
\section*{Introduction}

In the study of dynamical systems several authors considered the 
problem of determining which compact manifolds $M$ admit expansive systems. 
Let us recall that in the discrete-time setting 
$f\colon M\to M$ is an \emph{expansive homeomorphism} \cite{Utz} if 
there is $\delta>0$ such that $\dist(f^n(x),f^n(y))<\delta$ for all $n\in\Z$ implies $x=y$. 
In \cite{JU} it is proved that the circle does not admit expansive homeomorphisms.
In \cite{OR} it is shown that every orientable compact 
{surface} of positive genus admits expansive homeomorphisms, namely, 
a pseudo-Anosov {diffeomorphism}.
In \cites{Hir,Lew} it is
proved that the two-sphere does not admit expansive homeomorphisms. 
They also proved that every expansive surface 
homeomorphism is conjugate to a pseudo-Anosov {diffeomorphism}. 
This completes a global picture of expansive homeomorphisms of orientable compact surfaces.

{In higher dimensions} there is no result characterizing
which manifolds admit expansive homeomorphisms. 
Let us mention some advances in this direction. 
In \cite{Hi89} it is proved that expansive homeomorphisms of tori with the pseudo-orbit tracing property 
are conjugate to hyperbolic automorphisms.
In \cites{V93,V962} it is proved 
that an expansive 
homeomorphism of a compact three-dimensional manifold with 
a dense set of topologically hyperbolic periodic points 
is conjugate to a linear Anosov isomorphism on the torus. 
In \cite{ABP} this result was generalized for arbitrary dimension {assuming the existence of a codimension one periodic point}.
In \cite{Vi2002} it is proved {that 
expansive} $C^{1+\theta}$-diffeomorphisms 
on three-manifolds without wandering points are {conjugate to} Anosov diffeomorphisms on the torus.
The main difficulty, from our viewpoint, 
for classifying expansive homeomorphisms of three-manifolds is to understand the topology 
of local stable and unstable sets. 
To our best knowledge it is unknown whether the sphere
$\Stres$ admits an expansive homeomorphism. 

For the case of vector fields or flows 
the corresponding problems are considered.
According to Bowen and Walters \cite{BW} we say that 
$\phi\colon \R\times M\to M$ is an \emph{expansive flow} if for all $\epsilon>0$ there is $\delta>0$ such that 
if $\dist(\phi_{h(t)}(x),\phi_t(y))<\delta$ for all $t\in \R$ being $h\colon\R\to\R$ a parameterization, i.e., 
an increasing homeomorphism 
with $h(0)=0$, then $y=\phi_s(x)$ for some $s\in(-\epsilon,\epsilon)$.
An important fact about this definition is that it does not allow singular (or equilibrium) {points. 
Trivially}, every circle flow without singular points is expansive. 
It is known that no compact surface admits an expansive flow \cites{LG,Ar,CM} in the sense of \cite{BW}. 
In \cite{Pa} it is shown that if a compact three-manifold admits an expansive flow then its fundamental group 
has exponential growth. In particular, the three-sphere does not admit expansive flows in the sense of Bowen and Walters.

The expansiveness of flows with singular points was first investigated in \cite{K}. 
In this paper Komuro proved that the Lorenz attractor is $k^*$-expansive, a definition designed 
to allow singularities. 
According to \cite{K}, a flow is $k^*$-\emph{expansive} if for all $\epsilon>0$ there is $\delta>0$ 
such that if $\dist(\phi_{h(t)}(x),\phi_t(y))<\delta$ for all $t\in \R$ being $h$ a reparameterization,
then $\phi_{h(t_0)}(x)=\phi_{t_0+s}(y)$ for some $s\in(-\epsilon,\epsilon)$ and $t_0\in\R$.
In \cite{Ar} it is proved that a flow is $k^*$-expansive if and only if 
for all $\epsilon>0$ there is $\delta>0$ such that 
if $\dist(\phi_{h(t)}(x),\phi_t(y))<\delta$ for all $t\in \R$ and a reparameterization $h$, 
then there is $z\in M$ such that $x,y\in\phi_{[0,\epsilon]}(z)$ and $\diam(\phi_{[0,\epsilon]}(z))<\epsilon$.
It is easy to see that a circle flow is $k^*$-expansive if and only if it has a finite number of singularities. 
In \cite{Ar} it is shown that every $k^*$-expansive surface flow is obtained from surgery on the suspension 
of minimal interval exchange maps. 
It is also proved that a surface admits a $k^*$-expansive flow
if and only if it is a two-torus with $b$ boundary components, $h$ handles and $c$ cross-cups with $b + h + c > 0$. In
particular the two-torus and the two-sphere do not admit $k^*$-expansive flows.

From the definitions it is easy to see that every expansive flow in the sense of 
Bowen and Walters is $k^*$-expansive.
In order to obtain $k^*$-expansive flows with singular points on a manifold 
with dimension greater than 2 we can proceed as follows. 
Take $M$ admitting a Bowen-Walters expansive flow generated by a vector field $X$. 
Let $\rho\colon M\to \R$ be a non-negative smooth function vanishing only at $p\in M$. 
The vector field $\rho X$ generates a $k^*$-expansive flow with a zero-index singular point $p$. 
These kind of points are usually called \emph{fake singularities}. 
No published example of a $k^*$-expansive flow with hyperbolic singularities 
on a manifold of dimension greater than 2 is known to the author.
{This kind of expansiveness was deeply studied 
in relation with singular hyperbolic vector fields and three-dimensional attractors, 
as for example the one discovered by Lorenz, see \cite{AP}. 
Also, the concept of sectional-Anosov flow seems to be related with $k^*$-expansivity \cite{BM}.}

The purpose of the present paper is to show that the three-sphere $\Stres$ admits $k^*$-expansive flows. 
Let us sketch the construction while describing the contents of the article.
In Section \ref{PS} we will consider a triangular billiard in the hyperbolic disc, i.e., the curvature of the surface is -1 
and the boundary consists of three geodesic arcs. 
This dynamical system is related with the geodesic flow of a two-sphere with three punctures. 
The unit tangent bundle of this three-punctured sphere will be embedded in a closed three-manifold $M$. 
In Section \ref{secTopo} it is shown that $M$ is homeomorphic to $\Stres$. 
In Section \ref{SGF} we will show that a reparameterization of the geodesic flow of the three-punctured sphere 
can be extended to the whole $M$. 
In Section \ref{Exp} we show that this extended flow is $k^*$-expansive. 
{In Corollary \ref{corTrans} we deduce that the three-sphere admits a transitive $k^*$-expansive 
flow with a dense set of periodic orbits.}


\section{The phase space}
\label{PS}
In this section we will construct the phase space manifold $M$ of our $k^*$-expansive flow. 
It will be defined as a compactification of the unit tangent bundle of 
a three-punctured sphere {$\Sup^*$}. 
In Section \ref{secPhaseSpace} we will construct a smooth structure covering the punctures. 
In Section \ref{secExtPhaseSpace} we define a smooth manifold $M$ 
by adding three circles to {$T^1\Sup^*$}.
\subsection{A three-punctured sphere}
\label{secPhaseSpace}
\begin{figure}[htbp]
\begin{center}
   \includegraphics{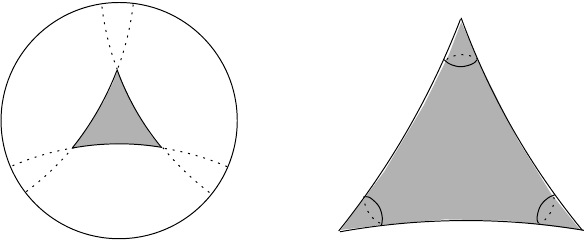}
    \caption{On the left, the hyperbolic triangle $\T\subset \D$. 
    On the right, the three-punctured sphere $\Sup$ associated to the triangle $\T$.}
   \label{figTdoble}
\end{center}
\end{figure}

Let $\D$ be the Poincaré disc with constant curvature -1.
Consider a triangle $\T\subset\D$ with geodesic sides as in Figure \ref{figTdoble}. 
Let $\varphi\colon \D\to \D$ be an isometry such that $\varphi(\T)\cap \T=\emptyset$. 
Let $\T_2=\T\cup \varphi(\T)$.
Define 
$\Sup=\T_2/\simeq,$
where $\simeq$ is the equivalence relation 
on $\T_2$ generated by $x\simeq \varphi(x)$ for all $x\in \partial \T$.
That is, we are gluing the boundaries of two disjoint copies of $\T$ as shown in Figure \ref{figTdoble}.
With the quotient topology, $\Sup$ is homeomorphic to the two-dimensional sphere $\Sdos$.
The equivalence class of a vertex of $\T$ is called as a \emph{conical singularity} of $\Sup$.
Denote by $\vertex$ the set of singular points of $\Sup$.
Since the boundary of $\T$ is made of geodesic arcs, we have 
that $\Sup^*=\Sup\setminus\vertex$ admits a 
natural smooth structure with a Riemannian metric of curvature -1. 
{This Riemannian metric induces a distance in $\Sup^*$ that extends to $\Sup$ and will be called
$\dist$.}
For $\radio>0$ 
define the ball and the reduced ball, respectively, as
\begin{equation}
\label{Bsigma}
\begin{array}{l}
B^*_\radio(q)=\{p\in \Sup:0<\dist(p,q)<\radio\},\\
B_\radio(q)=\{p\in \Sup:0\leq\dist(p,q)<\radio\}. 
\end{array}
\end{equation}
Fix $\radio>0$ so that 
for every singular point $\sigma\in \vertex$ 
the closure of 
$B_\radio(\sigma)$
is homeomorphic to a compact disc and $B_\radio(\sigma)\cap\vertex=\{\sigma\}$.

\begin{figure}[htbp]
\begin{center}
   \includegraphics{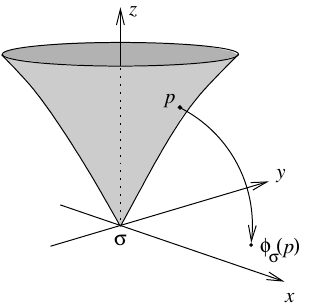}
    \caption{Geometry of the chart $\phi_\sigma$.}
   \label{Mapaphi}
\end{center}
\end{figure}

We will define polar coordinates around $\sigma$. 
Consider the map $r\colon B_\radio(\sigma)\to [0,\radio)$ given by 
$r(p)=\dist(p,\sigma)$.
Fix a geodesic $l$ starting at $\sigma$.
Denote by $\theta$ the angle at $\sigma$ of the triangle $\T$, i.e., $2\theta$ is the angle of the cone at $\sigma$. 
Let us denote as $\Suno_\theta$ the circle $\R/2\theta\Z$, for $\theta>0$, and as a special case 
 $\Suno=\R/2\pi\Z$.
\begin{figure}[htbp]
\begin{center}
   \includegraphics{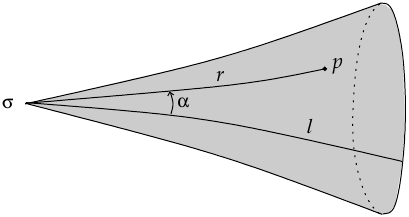}
    \caption{Polar coordinates near a singular point $\sigma$ in the surface $\Sup$.}
   \label{coord1}
\end{center}
\end{figure}
Define $\alpha\colon B^*_\radio(\sigma)\to\Suno_\theta$
such that $\alpha(p)$ is the angle 
between the geodesic segment from $\sigma$ to $p$ and $l$
as shown in Figure \ref{coord1}. 
Consider the isomorphism 
\begin{equation}
 \label{zetaSuno}
\zeta_\theta \colon\Suno_\theta\to\Suno \hbox{ such that } 
\zeta_\theta(\overline{x})=\overline{x\pi/\theta},
\end{equation}
for $x\in\R$, where the line over a point denotes the class in the corresponding quotient.
Consider the plane disc
\[
 \begin{array}{l}
  D_\radio=\{(x,y)\in\R^2:\sinh(\sqrt{x^2+y^2})\in [0,\radio)\}.
 \end{array}
\]
and define the map $\phi_\sigma\colon B_\radio(\sigma)\to D_\radio$
as
\begin{equation}\label{difeo}
\phi_\sigma(p)=
\sinh(r(p))(\cos(\zeta_\theta\circ\alpha(p)),\sin(\zeta_\theta\circ\alpha(p))) 
\end{equation}
if $p\neq \sigma$
and $\phi_\sigma(\sigma)=(0,0)$.
In Figure \ref{Mapaphi} we illustrate the map $\phi_\sigma$ as if $B_\radio(\sigma)$ were embedded in Euclidean $\R^3$. 
\footnote{At this point it may not be clear the reason why we use $\sinh(r)$ in (\ref{difeo}) instead of $r$, as would 
seem natural. This choice is done in order to simplify Section \ref{SGF}.}
Since $\phi_\sigma$ restricted to $B^*_\radio(\sigma)$ is a {diffeomorphism}, 
the charts $\{\phi_\sigma:\sigma\in\vertex\}$ extend the smooth atlas of $\Sup^*$ to $\Sup$.
This smooth structure of $\Sup$ and the polar coordinates around the singular points will be used in the following sections.

\subsection{The extended phase space}
\label{secExtPhaseSpace}
Define $M^*=T^1\Sup^*$, the unit tangent bundle of $\Sup^*$. 
Recall from the previous section that $\Sup^*$ is a three-punctured sphere.
In this section we will construct a closed three-manifold $M$ so that
$M^*\subset M$ 
and $M\setminus M^*$ is the disjoint union of three circles.

\begin{figure}[htbp]
\begin{center}
   \includegraphics{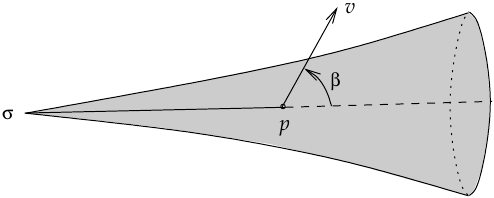}
    \caption{The map $\beta$.}
   \label{coord2}
\end{center}
\end{figure}

Denote by $\Pi\colon M^*\to \Sup^*$ the canonical projection.
Recall, from equation (\ref{Bsigma}), that $B^*_\radio(\sigma)$ is a reduced ball in $\Sup$ around the singularity $\sigma$. 
Define
$
U^*_\radio(\sigma)=\Pi^{-1}(B^*_\radio(\sigma))
$
and $\beta\colon U^*_\radio(\sigma)\to\Suno$
such that $\beta(v)$ is the angle between the tangent vector $v\in T^1_x\Sup$ 
and the geodesic from $\sigma$ to $x$ 
as shown in Figure \ref{coord2}. 
The direction of $\beta$ must be coherent with a fixed orientation of $\Sup$.
For each $\sigma\in\vertex$ consider the circle $\gamma_\sigma=\{\sigma\}\times \Suno$. 
Define $U_\radio(\sigma)=U^*_\radio(\sigma)\cup \gamma_\sigma$ and the map
\begin{equation}\label{cartasing}
\Phi_\sigma\colon U_\radio(\sigma)\to D_\radio\times\Suno 
\end{equation}
as
\[
 \left\{
 \begin{array}{ll}
\Phi_\sigma(v)=(\phi_\sigma(\Pi(v)),\beta(v)) & \hbox{ if } v\in U^*_\radio(\sigma),\\
\Phi_\sigma(\sigma,\beta)=(0,0,\beta)& \hbox{ if } (\sigma,\beta)\in\gamma_\sigma.
\end{array}
 \right.
\]
Recall that $\phi_\sigma$ was defined in equation (\ref{difeo}).
Consider the set
$$M=M^*\cup\bigcup_{\sigma\in\vertex}\gamma_\sigma.$$ 
The smooth structure of $M$ is defined via the maps $\Phi_\sigma$.
Therefore, the inclusion $M^*\to M$ is a diffeomorphism. 
Since $M$ can be covered by a finite number of charts with compact image, 
we have that $M$ is a compact three-dimensional manifold. 
Notice that by construction the manifold $M$ has no boundary.

\section{The topology of the phase space}
\label{secTopo}
We will study the topology of the manifold $M$. 
In Section \ref{secPrincCircBundle} we show that $M$ is a principal circle bundle. 
{In Section \ref{secSimpCon} we show that $M$ is homeomorphic to $\Stres$.}

\subsection{Principal circle bundle structure}
\label{secPrincCircBundle}

Consider the sphere $\Sup$ with singular set $\vertex$. 
On $\Sup$ consider the smooth structure given in Section \ref{secPhaseSpace}.
Recall that $M^*$ is the unit tangent bundle of $\Sup^*=\Sup\setminus\vertex$ and 
$\gamma_\sigma=\{\sigma\}\times \Suno$ for $\sigma\in\vertex$. 
The circle $\Suno$, as a Lie group, will be considered with additive notation $(\Suno,+)$.
In Section \ref{secExtPhaseSpace} we constructed 
a smooth structure for $M=M^*\cup_{\sigma\in\vertex}\gamma_\sigma$. 
Let us extend the canonical projection
$\Pi\colon M^*\to \Sup^*$  to
$\Pi\colon M\to \Sup$ 
as $\Pi(\sigma,\beta)=\sigma$ for $\sigma\in \vertex$.
We will show that we have a principal bundle structure 
\[
 \Suno\to M\to \Sup
\]
with structure group $\Suno$. 
\begin{prop}
The map $\Pi\colon M\to\Sup$ is a smooth submersion. 
\end{prop}

\begin{proof}

We already know that $\Pi$ restricted to $M^*$ is smooth.
Thus, let us consider a singular point $\sigma\in\vertex$ and 
the local charts $\phi_\sigma\colon B_\radio(\sigma)\to D_\radio$ and 
$\Phi_\sigma\colon U_\radio(\sigma)\to D_\radio\times \Suno$ given in (\ref{difeo}) and (\ref{cartasing}) respectively. 
If we define $P_1\colon D_\radio\times \Suno\to D_\radio$ as the projection on the first coordinate, 
we have the following commuting diagram:

\begin{center}
\begin{tikzpicture}
\matrix(m)[matrix of math nodes,
row sep=2.6em, column sep=2.8em,
text height=1.5ex, text depth=0.25ex]
{U_\radio(\sigma) & D_\radio\times \Suno\\
B_\radio(\sigma) &D_\radio \\};
\path[->,font=\scriptsize,>=angle 90]
(m-1-1) edge node[auto] {$\Phi_\sigma$} (m-1-2)
edge node[auto] {$\Pi$} (m-2-1)
(m-1-2) edge node[auto] {$P_1$} (m-2-2)
(m-2-1) edge node[auto] {$\phi_\sigma$} (m-2-2);
\end{tikzpicture}
\end{center}
Since $P_1$ is a smooth submersion the same is true for $\Pi$.
\end{proof}

For $p\in\Sup$ define $\gamma_p=\Pi^{-1}(p)$ as the \emph{fiber} at $p$.
We will define an action $\rt\colon\Suno\times M\to M$ 
of $\Suno$ on $M$. 
For $\varphi\in\Suno$ 
define $\rt_\varphi \colon T^1_x\Sup^*\to T^1_x\Sup^*$ 
as the rotation of angle $\varphi$. 
The direction of the rotation 
is determined by the orientation of $\Sup$ that we fixed in Section \ref{secExtPhaseSpace} in order to define the angle $\beta$.
For a point $(\sigma, \beta)\in\gamma_\sigma$ in a singular fiber define 
$\rt_\varphi(\sigma,\beta)=(\sigma,\varphi + \beta)$. 

\begin{obs}
 The circle action $\rt\colon\Suno\times M\to M$ is free, i.e. $\rt_\varphi(v)=v$ implies $\varphi=0$, smooth 
 and preserves the fibers of $\Pi \colon M\to \Sup$.
These remarks are obvious for $M^*$. 
The result near a singular fiber can be proved using local coordinates. 
Around a singularity we have that $$\beta(\rt_\varphi(v))=\beta(v)+\varphi.$$ 
In fact, this proves that 
\[
 \Suno\to M\to\Sup
\]
is a principal circle bundle with structure group $\Suno$.
\end{obs}

\subsection{The fundamental group of the phase space}
\label{secSimpCon}
{In this section we will show that $M$ is homeomorphic to the three-sphere.}
Consider the two-sphere $\Sup$ with its smooth structure defined above. 
For a smooth vector field $Y$ on $\Sup$ 
denote by $\sing(Y)$ the set of singular points (equilibrium) of $Y$. 
A singular point is a \emph{perfect center} of $Y$ if 
every trajectory of $Y$ near $\sigma$ is  
the boundary of a ball of $\dist$ centered at $p$. 
{In Figure \ref{coord2} we can see that a perfect center is characterized by $\beta=\pm\pi/2$.}

\begin{teo}
\label{piuno}
Each fiber of 
 $
  \Suno\to M\to\Sup
 $
bounds a two-dimensional disc. 
\end{teo}

\begin{proof}
Let $\sigma_1,\sigma_2,\sigma_3$ be the singularities of $\Sup$.
We will construct a disc whose boundary is the fiber at $\sigma_1$. 
Since all the fibers are homotopic this is sufficient.
Recall from equation (\ref{cartasing}) the local charts around the singular fibers 
$\Phi_\sigma\colon U_\radio(\sigma)\to D_\radio\times\Suno$.
Consider a smooth vector field $Y\colon\Sup^*\to M^*$ satisfying:
\begin{enumerate}
  \item $\Phi_{\sigma_1}\circ Y(x)=(\phi_{\sigma_1}(x), -\zeta_\theta(\alpha(x)))$ for $x$ close to $\sigma_1$, 
  recall the map $\alpha$ from the polar coordinates illustrated in Figure \ref{coord1}. The map $\zeta_\theta\colon \Suno_\theta\to\Suno$ 
  is the isomorphism defined in Equation (\ref{zetaSuno}). 
  {The last coordinate in $D_\radio\times\Suno$ represents the angle $\beta$, therefore, we are requiring that 
  $\beta(Y(x))=-\zeta_\theta(\alpha(x))$.}
  \item $\Phi_{\sigma_2}\circ Y(x)=(\phi_{\sigma_2}(x), \pi/2)$ for $x$ close to $\sigma_2$. 
  {In this case we have $\beta(Y(x))=\pi/2$.}
  \item $\Phi_{\sigma_3}\circ Y(x)=(\phi_{\sigma_3}(x),-\pi/2)$ for $x$ close to $\sigma_3$, {which implies $\beta(Y(x))=-\pi/2$ near $\sigma_3$.}
\end{enumerate}
The solutions of such vector field are illustrated in Figure \ref{campoY}. 
\begin{figure}[htbp]
\begin{center}   
    \includegraphics{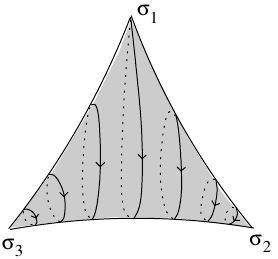}
    \caption{Flow lines of the vector field $Y$.}
    \label{campoY}
\end{center}
\end{figure}
Note that conditions 2 and 3 means that $\sigma_2$ and $\sigma_3$ are perfect centers. 
{As we said, the coordinate $\beta$ is constant near these singularities.
Therefore $Y$ can be extended to $Y'\colon\Sup\setminus\{\sigma_1\}\to M\setminus \gamma_{\sigma_1}$. }
We have that $Y'\colon\Sup\setminus\{\sigma_1\}\to M$ is a 
homeomorphism onto its image.
Notice that $\Sup\setminus\{\sigma_1\}$ is a one-punctured sphere and therefore it is homeomorphic to a an open disc. 
Condition 1 implies that the boundary of the disc $Y'(\Sup\setminus\{\sigma_1\})$ is the fiber at $\sigma_1$.
\end{proof}

\begin{cor}
The fundamental group $\pi_1(M)$ is trivial and the principal bundle structure $ \Suno\to M\to\Sup
$ is equivalent to the Hopf bundle. 
In particular, $M$ is homeomorphic to $\Stres$.
\end{cor}

\begin{proof}
{ It is known, see for example \cite{Naber}, that every principal circle bundle
 $\Suno\to N\to \Sdos$, 
 with homotopically trivial fibers
 is equivalent with the Hopf bundle. 
Thus, the result follows by 
Theorem \ref{piuno}.}
\end{proof}

\section{The extension of the flow}
\label{SGF}
In this section we will show that the geodesic flow of $\Sup^*$, defined on $M^*=T^1\Sup^*$, can be reparameterized and extended to $M$.
Let $X\colon M^*\to T M^*$ be the velocity field of the geodesic flow of $\Sup^*$.
Consider $\epsilon>0$ and 
a positive smooth function $\rho\colon M^*\to\R$ 
such that if $v\in M^*$, $\sigma\in\vertex$ and $\dist(\Pi(v),\sigma)<\epsilon$ then 
$\rho(v)=\sinh(\dist(\Pi(v),\sigma))$.
Let us define a time change of $X$ as $Y=\rho X$.
The main result of this section is the following:

\begin{teo}
\label{teoExtBilFlo}
The vector field $Y$ on $M^*$ has a smooth extension to $M$ that will be called $Z$.
\end{teo}

{The following remarks may help in the understanding of the geometric ideas behind the proof of this theorem.}

\begin{obs}
{We will prove that on each singular fiber there are two hyperbolic singularities of $Z$.
As we can see in Figure \ref{coord2}, the set of points of $M$ converging to the vertex are determined by the condition 
$\beta=\pi$. These points form the stable manifold of one of the singularities of $Z$ that will appear in $M$. 
The condition $\beta=0$ corresponds to the unstable manifold of the other singularity.
The complement of the singularities in the singular fiber consists of two trajectories connecting them.}
\end{obs}

In order to prove {Theorem \ref{teoExtBilFlo}} we give a preliminary result.
Recall the coordinates $(r,\alpha,\beta)$ 
introduced in Section \ref{PS}.

\begin{lemma}
The velocity field $X$ of the geodesic flow in the coordinates $(r,\alpha,\beta)$ 
 near a singular point $\sigma\in\Sup$, {denoted as $X_\sigma$,}  
is given by
$$X_\sigma(r,\alpha,\beta)=  
\left(\cos\beta,\frac1{\sinh r}\sin\beta ,-\frac{\cosh r}{\sinh r}\sin\beta\right).$$
\end{lemma}

\begin{proof}
Denote by $\phi$ the geodesic flow of $\Sup$. 
Consider $(r_0,\alpha_0,\beta_0)$ an initial condition and denote $\phi_t(r_0,\alpha_0,\beta_0)=(r_t,\alpha_t,\beta_t)$. 
\begin{figure}[htbp]
\begin{center}   
    \includegraphics{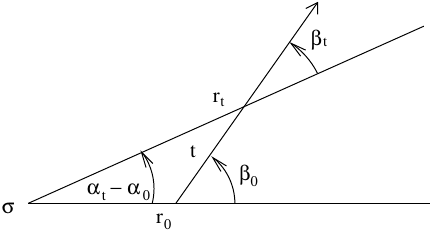}
    \caption{Geodesic flow near $\sigma$.}
    \label{billflow}
\end{center}
\end{figure}
Applying the laws of sines and cosines of the hyperbolic trigonometry 
on the triangle of Figure 
\ref{billflow} we conclude the following equations:
\[
\left\{
 \begin{array}{l}
  \sin \beta_t \sinh t  =\sin ({\alpha_t-\alpha_0}) \sinh r_0,\\
  \sin \beta_t \sinh r_t=\sin (\pi-\beta_0) \sinh r_0,\\
  \cosh r_t=\cosh r_0\cosh t- \sinh r_0\sinh t\cos (\pi-\beta_0).
 \end{array}
\right.
\]
Taking derivatives at $t=0$ we obtain:
\[
\left\{
 \begin{array}{l}
    \dot r=\cos\beta,\\
     {\dot \alpha}=\frac1{\sinh r}\sin\beta,\\
    \dot \beta=-\frac{\cosh r}{\sinh r}\sin\beta.
 \end{array}
\right.
\]
Therefore, the vector field 
$$X_\sigma(r,\alpha,\beta)=\left(\cos\beta,\frac1{\sinh r}\sin\beta ,-\frac{\cosh r}{\sinh r}\sin\beta\right)$$
is the velocity field of the geodesic flow in the coordinates $(r,\alpha,\beta)$.
\end{proof}

\begin{proof}[Proof of Theorem \ref{teoExtBilFlo}]
Fix a singular point $\sigma\in\vertex$. 
Consider {$X_\sigma$} the expression of 
$X$ in coordinates $(r,\alpha,\beta)$ near the fiber $\gamma_\sigma$.
Since $Y_\sigma=\rho X_\sigma$, we have that
$$Y_\sigma(r,\alpha,\beta)=(\sinh (r) \cos(\beta),\sin(\beta),-{\cosh(r)}\sin(\beta))$$
on local charts. 
Recall that the smooth structure of $M$ was defined with the charts $\Phi_\sigma$ in 
equation (\ref{cartasing}).
Consider the diffeomorphism
\[
 \Phi'_\sigma\colon (0,\epsilon)\times \Suno_\theta\times \Suno\to {D^*_\tau}\times\Suno
\]
defined as 
\begin{equation}
 \label{Phi'}
\Phi'_\sigma(r,\alpha,\beta)=(\sinh(r)\cos(\zeta_\theta(\alpha)),\sinh(r)\sin(\zeta_\theta(\alpha)),\beta). 
\end{equation}
It is the diffeomorphism $\Phi_\sigma$ in coordinates $(r,\alpha,\beta)$.
Define $Z_\sigma=\dif\Phi'_\sigma (Y_\sigma)$, a vector field in ${D^*_\tau}\times\Suno$. 
Let $(x,y,z)=\Phi'_\sigma(r,\alpha,\beta)$. 
In the coordinates $(x,y,z)$ the expression of {$Z_\sigma=Z_\sigma(x,y,z)$} is:
\begin{equation}\label{fieldflow}
 {Z_\sigma=
\cos z\sqrt{1+x^2+y^2}(x, y,0) - \sin z(y\pi/\theta ,-x\pi/\theta,\sqrt{1+x^2+y^2}}).
\end{equation}
This formula is proved in Lemma \ref{lemaAux} below.
We define $Z_\sigma$ in the whole 
$D_\radio\times\Suno$ by this formula, obtaining a smooth vector field.
Proceeding in the same way on each $\sigma\in\vertex$  
we can smoothly extend $Y$ to the whole of $M$.
\end{proof}

\begin{lemma}
\label{lemaAux}
 The expression of the vector field $Z_\sigma$ in coordinates $(x,y,z)$ is given by equation (\ref{fieldflow}).
\end{lemma}
\begin{proof}
Recall that $\zeta_\theta\colon\Suno_\theta\to\Suno$ was 
defined as the quotient map of $x\mapsto x\pi/\theta$. 
Then, by equation (\ref{Phi'}) we have that:
\[
\dif\Phi'_\sigma|_{(r,\alpha,\beta)}=
\left[
\begin{array}{ccc}
  \cosh(r)\cos(\zeta_\theta(\alpha)) & -\frac\pi\theta \sinh(r)\sin(\zeta_\theta(\alpha))& 0\\
  \cosh(r)\sin(\zeta_\theta(\alpha)) &  \frac\pi\theta \sinh(r)\cos(\zeta_\theta(\alpha))& 0\\
  0 & 0 & 1
\end{array}
\right]
\]
Recall that $Z_\sigma=\dif\Phi'_\sigma (Y_\sigma)$. Therefore
\[
  Z_\sigma= 
\left[
\begin{array}{c}
  \cosh(r)\cos(\zeta_\theta(\alpha))\sinh (r) \cos(\beta) -  \frac\pi\theta \sinh(r)\sin(\zeta_\theta(\alpha))\sin(\beta)\\
  \cosh(r)\sin(\zeta_\theta(\alpha))\sinh (r) \cos(\beta) +  \frac\pi\theta \sinh(r)\cos(\zeta_\theta(\alpha))\sin(\beta)\\
  -{\cosh(r)}\sin(\beta)
\end{array}
\right]
 \]
 Since $x=\sinh(r)\cos(\zeta_\theta(\alpha))$ and $y=\sinh(r)\sin(\zeta_\theta(\alpha))$ 
 we have that $\sqrt{x^2+y^2}=\sinh(r)$. Recalling that $\cosh(\sinh^{-1}(a))=\sqrt{1+a^2}$ for all $a\in\R$ we have that 
 $\cosh(r)=\sqrt{1+x^2+y^2}$. With these equations and $z=\beta$ 
 the formula (\ref{fieldflow}) is obtained.
\end{proof}

\begin{prop}
\label{propHypSing} 
The extended vector field $Y$ has two hyperbolic singularities on each singular fiber.
\end{prop}

\begin{proof}
 By construction, the equilibrium points of the flow are in 
the singular fibers $\gamma_\sigma$ for $\sigma\in \vertex$. So we consider the coordinates $(x,y,z)$ around $\sigma$ 
and the expression of the vector field given by (\ref{fieldflow}). Since 
$
Z_\sigma(0,0,z)=\left(0,0,-\sin z\right),
$
on each $\gamma_\sigma$ there are two equilibrium points: $p_1=(0,0,0)$ and $p_2=(0,0,\pi)$. 
To simplify the notation define $f(x,y)=\sqrt{1+x^2+y^2}$.
The linear part of $Z_\sigma$ in coordinates $(x,y,z)$ is
\newcommand{\aux}{f(x,y)}

\[
 \left[
 \begin{array}{c|c|c}
  \frac{1+2x^2+y^2}{\aux}\cos z              & \frac{xy}{\aux}\cos z-\frac\pi\theta\sin z & -x\aux\sin z-\frac\pi\theta y\cos z\\
  \frac{xy}{\aux}\cos z+\frac\pi\theta\sin z & \frac{1+x^2+2y^2}{\aux}\cos z              & -y\aux\sin z+\frac\pi\theta x\cos z\\
  {-\frac{x\sin z}{\aux}}    & {-\frac{y\sin z}{\aux}}    & {-\aux\cos z}
 \end{array}
 \right].
\]
Then 
\[
 d_{p_1}Z_\sigma=
 \left[
 \begin{array}{ccc}
 1 & 0 & 0\\
 0 & 1 & 0\\
 0 & 0 & -1
 \end{array}
\right]
\]
and
\[
d_{p_2}Z_\sigma=
 \left[
 \begin{array}{ccc}
 -1 &  0 & 0\\
 0  & -1 & 0\\
 0  &  0 & 1
 \end{array}
\right].
\]
Therefore $p_1$ and $p_2$ are hyperbolic singularities.
\end{proof}

\section{Dynamical properties of the flow}
\label{Exp}

In this section we will prove that the vector field $Z$ of Theorem \ref{teoExtBilFlo} defines a flow $\phi$ in $M$ that is $k^*$-expansive.
Since we have proved that $M$ is homeomorphic to $\Stres$ the result of the paper will be proved. 
Recall that a flow is $k^*$-\emph{expansive} if for all $\epsilon>0$ there exists an \emph{expansive constant} $\delta>0$ such that if 
$\dist(\phi_{h(t)}(x),\phi_t(y))<\delta$ for all $t\in\R$, being $h\colon\R\to\R$ an increasing homeomorphism with $h(0)=0$, 
then $x$ and $y$ are contained in an orbit segment of diameter
less than $\epsilon$.

\begin{teo}
\label{teoflat}
The flow $\phi$ associated to the vector field $Z$ is $k^*$-expansive in the {three-sphere} $M$.
\end{teo}

\begin{proof}
Given $\epsilon>0$ we will construct an expansive constant $\delta>0$ for $\phi$. 
First consider a regular point $x\in M$ for the flow. 
Take a local cross section $\sect_x$ and 
a flow box $B_x=\phi_{(-t_x,t_x)}\sect_x$ such that the diameter of every orbit segment contained in $B_x$ is smaller than $\epsilon$. 
{By Proposition \ref{propHypSing} we know that the singularities are hyperbolic.
Then, around each singular point $p$ of the flow we can consider 
an adapted neighborhood $B_p$, as in Figure \ref{figentornosing}.
Again, we assume that} the diameter of every orbit segment contained in $B_x$ is smaller than $\epsilon$. 
\begin{figure}[htbp]
\begin{center}
   \includegraphics{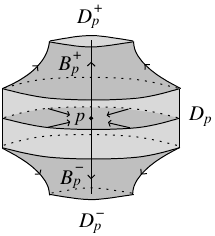}
    \caption{Adapted neighborhood $B_p$ of a singularity with two-dimensional stable manifold.}
   \label{figentornosing}
\end{center}
\end{figure}
For the singular point $p$ let $r=$ $s$ or $u$ be such that the stable (or unstable) manifold 
$W^r(p )$ has dimension 2. 
Define $D_p $ as the connected component of $W^r(p )\cap B_p $ that contains $p $. Let $D_p ^+$ and $D_p ^-$ be the local cross sections, contained in the boundary of 
$B_p $, given in Figure \ref{figentornosing}. 
Let $B_p ^+$ and $B_p ^-$ be the connected components of $B_p \setminus D_p $ such that 
$D_p ^+\subset\partial B_p ^+$ and $D_p ^-\subset\partial B_p ^-$, as in Figure \ref{figentornosing}. Define 
$$\delta_p =\min \{\dist(D_p ^+,B_p ^-),\dist(D_p ^-,B_p ^+)\}.$$
Since $M$ is compact we can take a finite covering of $M$ of the form 
$${\mathcal{ B}}=\{B_{x_1}, \dots,B_{x_n},B_{p_1}, \dots ,B_{p_6}\},$$
being $x_1,\dots,x_n$ some regular points and $p_1,\dots,p_6$ the singular points of the flow. 
Recall that there are three singular fibers in $M$, each one containing two singular points of the flow, giving 6 singularities.
Consider a positive $\delta$ such that $\delta<\min\{\delta_{p_i} :i=1,\dots,6\}$. 
Also assume that if $x,y\in M$ and 
$\dist(x,y)<\delta$ then there is $B\in\mathcal{B}$ containing $x$ and $y$.

We will show that $\delta$ is an expansive constant for the flow. 
By contradiction assume that $x,y\in M$ are not contained in an orbit segment of diameter 
$\epsilon$ and that there is a reparameterization $h$ such that $\dist(\phi_tx,\phi_{h(t)}y)<\delta$ for all $t\in\R$.
Recall that the flow was obtained as an extension of the geodesic flow of a three-punctured surface. 
Also, this surface was constructed by gluing two copies of the triangle 
$\T$ in the hyperbolic disc. 
Therefore we can view the dynamics of $\phi$ as the billiard flow of $\T$. 
By our choice of $\delta$ the billiard trajectories of $x$ and $y$ must have the same itinerary, 
i.e., the sequences of sides of $\T$ that they hit must coincide.
But this is impossible because the billiard surface has negative curvature and $x,y$ are in different local orbits.
\end{proof}

{Further properties of the flow can be deduced from the theory of billiards and geodesic flows.
We say that $\phi$ is a \emph{transitive} flow if there is $x\in M$ such that $\{\phi_t(x):t\in\R\}$ is dense in $M$.}

\begin{cor}
\label{corTrans}
  {The three-sphere admits a transitive 
  $k^*$-expansive flow with a dense set of periodic orbits.}
\end{cor}

\begin{proof}
In \cite{GSG} it shown that the billiard map of a triangle on the hyperbolic disk 
has non-vanishing Lyapunov exponents. 
Then, applying \cite{KS}*{Theorem 13.2} we have that 
the set of periodic orbits is dense in the phase space of the billiard map. 
This easily gives a dense set of periodic orbits in the three-sphere $M$ for the $k^*$-expansive 
flow $\phi$ of the previous theorem. 

In order to obtain an example presenting transitivity we will consider a special triangle. 
In the hyperbolic disk we can consider an equilateral triangle with angles $\pi/4$. 
A circular unfolding around a vertex gives us a regular octagon as in Figure \ref{figOctagon}.
\begin{figure}[h!]
\begin{center}
  \includegraphics[scale=.6]{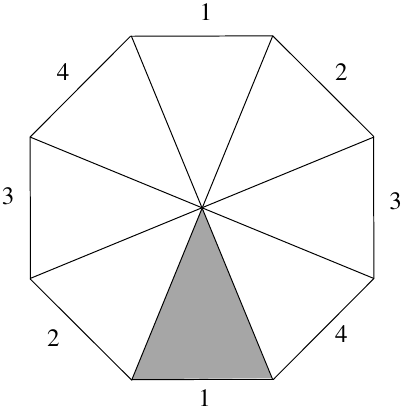}
  \caption{The unfolding of a triangle around a vertex. Identifying opposite sides a surface of genus two is obtained.}
  \label{figOctagon}
\end{center}
\end{figure}
If we identify opposite sides of this octagon we obtain a smooth 
surface of genus two with constant negative curvature.
The transitivity of its geodesic flow is a well known property, 
see for example \cite{BeMa}*{Corollary 3.8}. 
There, the ergodicity is proved for a measure that is positive on open sets, 
which easily implies the transitivity.
Now we must note that the set of points converging to a singular point in $M$ has vanishing Lebesgue measure. 
Therefore, there are dense trajectories of the geodesic flow of the genus two surface 
that projects into dense trajectories in $M$. 
This proves the transitivity of $\phi$ in the three-sphere $M$ obtained from the equilateral triangle of angles $\pi/4$.
\end{proof}

{Further properties of the flows obtained in this paper should be explored 
from the viewpoint of the ergodic theory as well as topological dynamics.
Also, it would be interesting to know which kind of knots appear as periodic orbits.
Let us finally indicate some possible extensions of our results.}

\begin{obs}
 {If instead of starting with a triangle we consider a polygon in the hyperbolic disc, we will obtain a $k^*$-expansive flow 
 but the ambient manifold will have non-trivial fundamental group. The reader can check this with the techniques of the paper.}
\end{obs}

\begin{obs}
 {If instead of considering a triangle with negative 
 curvature we start with a flat (Euclidean) triangle then 
 we can prove that the $k^*$-expansivity of the flow is equivalent 
 with the non-existence of periodic orbits in the triangular billiard. 
 This essentially follows by the results in \cite{GKT}.
 Let us mention that it is not known whether every flat triangular billiard has a periodic orbit or not.
 Hopefully, the ideas in the present article may help in the study of triangular billiards.}
\end{obs}

\end{document}